# WeBWorK log files as a rich source of data on student homework behaviours

Alain Prat[a] and Warren J. Code[b]

[a]Department of Mathematics, University of British Columbia, Vancouver, BC, Canada;
[b]Science Centre for Learning and Teaching (Skylight), University of British Columbia, Vancouver, BC, Canada

**ABSTRACT**
The online homework system WeBWorK has been successfully used at several hundred colleges and universities. Despite its popularity, the WeBWorK system does not provide detailed metrics of student performance to instructors. In this article, we illustrate how an analysis of the log files of the WeBWorK system can provide information such as the amount of time students spend on WeBWorK assignments and how long they persist on problems. We estimate the time spent on an assignment by combining log file events into sessions of student activity. The validity of this method is confirmed by cross referencing with another time estimate obtained from a learning management system. As an application of these performance metrics, we contrast the behaviour of students with WeBWorK scores less than 50% with the remainder of the class in a first year Calculus course. This reveals that on average, the students who fail their homework start their homework later, have shorter activity sessions, and are less persistent when solving problems. We conclude by discussing the implications of WeBWorK analytics for instructional practices and for the future of learning analytics in undergraduate mathematics education.

## 1. Introduction

Homework is a well-established practice in post-secondary mathematics courses, with a trend in recent decades to augment – or, in some cases, replace – paper-based homework with online homework systems, particularly in larger course offerings such as those in the Calculus sequence [1]. One of the most popular online systems, particularly in the United States, is WeBWorK; indeed, the Mathematical Association of American (MAA) website [2] states that 'WeBWorK is used successfully at over 700 colleges and universities from large research institutions to small teaching colleges. WeBWorK has been developed and maintained by mathematicians since 1994 always with the goal of providing the mathematical community with the most robust, flexible, and mathematically capable online homework system possible.' WeBWorK is free and open source, with an extensive library of problems that may be further customized as needed. Despite this appeal and popularity, and perhaps owing to its age, the WeBWorK system does not provide instructors with detailed information regarding students' interaction with the system [3]. In what follows, we will demonstrate how pedagogically useful information can be obtained from data already generated by this popular online homework system.

---

CONTACT Warren Code. Email: warcode@science.ubc.ca

Early WeBWorK studies were largely focused on the user experience of students and how it affected their overall approach to homework, with a goal of establishing an online homework system which would be at least as effective (if not more) compared to traditional paper-based homework in terms of student acceptance and student learning [4]. These studies included student surveys of perceived value, course performance, and other aspects of student experience and homework behaviour (for example, the likelihood of skipping difficult problems as related to course grade [5]), and this is still being revisited as the use of online homework has continued to expand in the teaching of undergraduate mathematics in WeBWorK [6,7] and other systems [8,9]. Over the years, experimental features to achieve new pedagogical aims have been evaluated, for example a collaborative digital whiteboard embedded in questions [10] or an enhanced system of hints for providing feedback [11]; these efforts have generally collected a new type of data originating from that same new feature, such as capturing mouse behaviour in the whiteboard example [10].

WeBWorK's web interface for instructors offers some descriptive statistics of student activity, including completion rates and an index of question difficulty based on attempts, as well as the ability to see past responses on a question-by-question basis for individual students. Here, we present a method for calculating additional metrics related to WeBWorK usage, using data that is routinely collected by the system. These metrics cannot currently be obtained from either the menus within WeBWorK or the scoring files. Rather, they are obtained by analyzing the pair of system log files which record logins and every answer submission submitted by every student. Using this log data, we demonstrate how it is possible to estimate the time students have spent on their WeBWorK assignments, as well as other timing metrics such as how long students persist on a single problem. These additional metrics have an untapped pedagogical potential within the context of WeBWorK. For example, contrasting the average time students spend on the different assignments in a course can allow an instructor to calibrate the relative difficulties of their assignments the following year. As another example, instructors can use the timing of answer submissions stored in the log files to determine how many students are starting their homework well before the due date, how many are starting at the last minute, and even how many are working on their WeBWorK assignment during classroom time that is meant to be devoted to other activities.

We were able to find examples of internal reports from institutions (including our own) that looked at student data from WeBWorK, in some cases self-published on a department website (e.g., [12]). However, these are primarily instructors reporting on instructional goals for their courses and seldom proceed to further publication. Data on completion and number of student attempts, available through the web interface for instructors, has been used as a way of analyzing problem difficulty [13]. It has also been suggested that this same data be used as part of a 'dashboard' (visualization tool) to help instructors monitor student success in the course's online homework [3]. This dashboard tool would augment what is currently available to instructors through the 'Statistics' tool that is part of the WeBWorK web interface.

The response log data we mention above contains the text of every student submission, and this has been used to analyze submission patterns to develop a 'Student Response Model' [5] that uses a team of human raters to classify response quality, for example, separating meaningfully new attempts to a problem from simply resubmitting the same or an equivalent mathematical expression, or gaming the response system by entering an obvious series of guesses. The response log has also been used to study response correctness via Item Response Theory and to produce a dashboard of student accomplishment and item difficulty [14]. However, to our knowledge, no published studies have used the information in the logged data which pertains to timing of student work in WeBWorK, apart from a simple use in estimating time spent on homework assignments as part of a larger study comparing teaching methods [15].

Similar to other research in online homework, we find that the students with the lowest assignment scores start their online mathematics homework later [16–18] and do not persist as long before giving up on a problem [5,18]. We also find that these low-score students have shorter activity sessions but attempt more problems per hour [18]. To our knowledge, these types of results have arisen from student surveys and interviews or log data in another online homework system, but not using data directly from the WeBWorK system.

Our work can be viewed in the context of the broader field of learning analytics in higher education, a collection of practices related to the explosion of data following the growth of online courses and technology enhanced learning [19]. While many studies have examined the potential of online interaction data to examine student time spent and engagement with online learning systems in general, for mathematics homework there is still substantial value to a purpose-built system (often connected to an LMS) and so it is natural to consider the data generated in these systems [3].

This article is organized as follows. In Section 2, we present the method that is used to extract student performance metrics from the WeBWorK log files. In Section 3, we illustrate how these WeBWorK metrics can be applied to contrast the behaviour of students with WeBWorK scores below 50% with the behavior of students in the remainder of the class. Finally, we discuss some of the implications of our work as well as future research directions in Section 4.

## 2. Methods

### 2.1. WeBWorK data sources

Currently, the most readily available data to instructors in WeBWorK is obtained in the 'Statistics' function of the Instructor Tools menu. For a given WeBWorK assignment, this tool provides instructors with the percentage of students who answered a problem correctly, among those who attempted that problem. Percentiles for overall assignment scores are also provided, as well as the average number of attempts per problem and selected percentile cutoffs for the number of attempts for each problem. This data can be used to compare the relative difficulty of problems within an assignment and potentially improve the assignment by addressing issues related to difficult problems. If an instructor desires additional information about the number of attempts and scores for each problem on an assignment, these can in principle be obtained from the Scoring Tools function of the Instructor Tools, which can generate a comma-separated values (csv) file containing the counts of correct and incorrect attempts for a each problem on an assignment. This data can provide instructors with additional information not contained in the Statistics menu, such as the percentage of students that did not attempt a problem, or the students in the class who required the most number of answer attempts to complete an assignment. However, in order to

```
[Fri Dec 02 23:01:13 2016] |AMFFP5X4I202|Assignment_12|10|0 1480748473    (sqrt(3)/2)+pi/12
[Fri Dec 02 23:01:19 2016] |76ARTLFSBF01|Assignment_12|24|1 1480748479    (120^2) / (32*2 )
[Fri Dec 02 23:01:34 2016] |KBGURC1AHF18|Assignment_12|9|0 1480748494     9*(-9^(1/3))/(1+-9^(1/3))
[Fri Dec 02 23:01:40 2016] |JT18Z8YBV504|Assignment_12|13|0 1480748499    ((h^2d-hd^2)/(h-d))^(1/2)
[Fri Dec 02 23:02:24 2016] |Q7DHCFXG8O09|Assignment_12|10|0 1480748544    2.6678
[Fri Dec 02 23:02:40 2016] |8MKDZ27AFT05|Assignment_12|26|0 1480748560    -18
[Fri Dec 02 23:04:00 2016] |DYRXI8W6ZC16|Assignment_12|5|00 1480748640    7.006   1/2(10-(1+pi/2)(20/(4+pi)))
[Fri Dec 02 23:04:39 2016] |CL9JMXD1PK09|Assignment_12|8|110 1480748679   3/2s^2csc^2(t)-(3s^2 (sqrt3)/2)cot(t)*csc(t)  cos^-1(sqrt3/3)
```

(a) Sample of lines from an answer log file.

```
[Wed Oct 26 13:47:33 2016] LOGIN OK user_id=6834XIFTZS03 login_type=normal credential_source=LTI host=123.456.789.9 port=40001
UA=Mozilla/5.0 (Macintosh; Intel Mac OS X 10_12) AppleWebKit/602.1.50 (KHTML, like Gecko) Version/10.0 Safari/602.1.50
[Wed Oct 26 13:48:32 2016] AUTH WWDB: password rejected, deferring to site_checkPassword user_id=1DWCV8NALJ04 login_type=normal
credential_source=params 123.456.789.9 port=40001 UA=Mozilla/5.0 (Windows NT 10.0; Win64; x64) AppleWebKit/537.36 (KHTML, like Gecko)
Chrome/53.0.2785.143 Safari/537.36
[Wed Oct 26 13:48:32 2016] LOGIN OK user_id=1EWCV9NALJ04 login_type=normal credential_source=params 123.456.789.9 port=40001
UA=Mozilla/5.0 (Windows NT 10.0; Win64; x64) AppleWebKit/537.36 (KHTML, like Gecko) Chrome/53.0.2785.143 Safari/537.36
```

(b) Sample of lines from a login file.

**Figure 1.** Examples of raw data obtainable from WeBWorK log files.

obtain this information, the raw comma-separated values file has to be processed by the instructor using a spreadsheet or other software.

The WeBWorK system records every correct and incorrect answer submission in a plain text system log file. The information in this file is used when an instructor uses the 'Show Past Answers' feature of WeBWorK, which allows the instructor to examine the record of an individual student's answers for a problem. The raw data can be obtained by the instructor navigating to 'Instructor Tools' – 'File Manager', locating the directory '/logs' in the home directory, and downloading the file 'answer log'. An example of several lines from such a log file is shown in Figure 1a. The columns are separated by a 'pipe' (vertical line) symbol, and provide the following data: date, user ID, assignment name, problem number, correctness of answer submission, time since Jan 1 1970 (in seconds), text of answer submitted. In the present work, we will only be using information from the first five columns, though as noted above the full file is also a source of data for analysis of answer patterns [5,14].

A second log file in WeBWorK records when a student successfully or unsuccessfully logs in, as well as when their sessions time out due to inactivity. An example of several lines from this log file is shown in Figure 1b. This file is found in the same directory as the answer log file mentioned above, and the file is named 'login.log'.

We have looked online to determine the awareness about and use of this information. As part of its support for WeBWorK use and development, the MAA hosts a web forum to serve as 'a place for users and administrators to collaborate and exchange expertise about the WeBWorK system' [20]. From this forum, one can learn that these files are known to system administrators of WeBWorK servers, and offer information to help diagnose problems with student access to their local WeBWorK site [21]. Log files have also been used by some instructors to recover lost information about student answer submissions, or to address technical issues related to answer submissions [22]. Beyond these administrative uses, a few examples exist of instructors interested in the type of pedagogically useful information that can be extracted from the log files [21,22]. This includes a 2009 posting looking for a script to estimate time spent on WeBWorK using the login.log file (Jan 2009 posting by R Gompa to WeBWorK Main Forum; unreferenced), and a 2018 posting asks if information obtained from the answer log file can be used to determine how students are reviewing WeBWorK assignments after their due date (Sep 2018 posting by J Trussell to WeBWorK Main Forum; unreferenced); both are unanswered as of this writing.

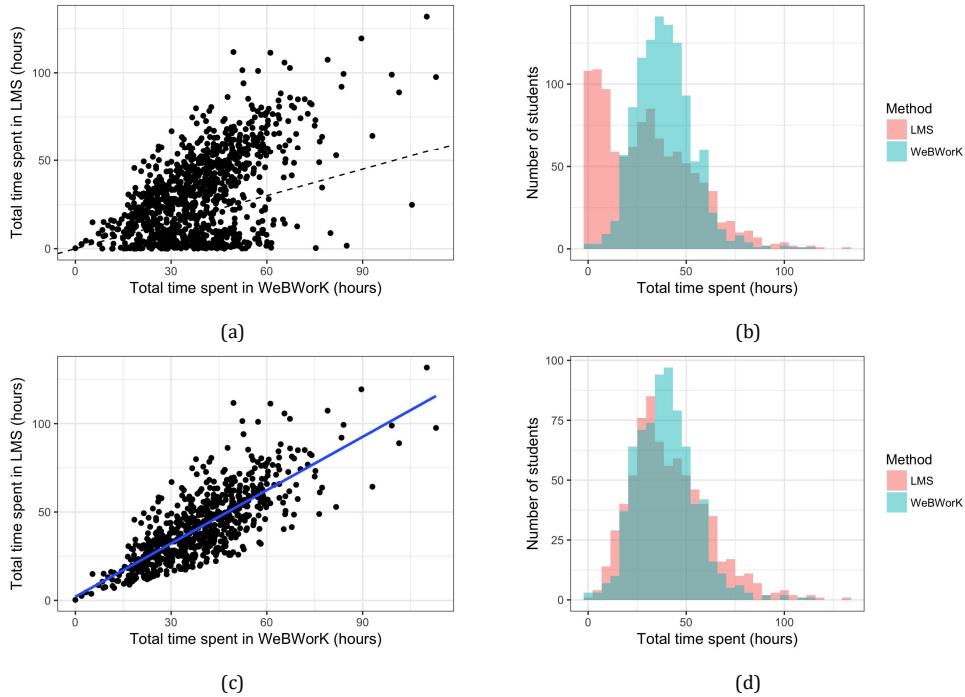

**Figure 2.** (a) Scatterplot showing the time spent in the LMS vs the estimated time spent in WeBWorK totalled over the whole course. Each point represents a student. Because it is possible to access WeBWorK without going through the LMS, there are a number of students with much less time in the LMS than in WeBWorK. (b) A comparison of the distributions for the time spent in the LMS and the time spent in WeBWorK. We see a bimodal pattern due to some students not going through the LMS to access WeBWorK. (c) The same scatterplot as in (a), but with the students whose WeBWorK time exceeds the LMS time by twice as much removed. (d) The distribution of LMS and WeBWorK times corresponding to the scatterplot in
(c).

## 2.2. Time on task estimation

By combining the data from the two log files, one can estimate the time a student spends working on a particular assignment. The algorithm we use for calculating this estimate is as follows: all logins and all answer submissions related to that assignment are considered events of activity. An activity session is defined by a sequence of such events with no gap between them lasting more than a predetermined inactivity threshold; if two events are separated in time by more than the inactivity threshold, we consider the first event to be the end of one activity session and the next event to be the beginning of the next activity session. A session length is the total time between the first and last events for that session, and the total time spent on an assignment is then the sum of the length of all the sessions associated with that assignment. We note that this algorithm is one of several possible methods of time-on-task estimation, and recent research has revealed that in some cases the choice of method can have profound effect on research findings [23].

Our threshold of inactivity was chosen by cross referencing the total time spent in WeBWorK, as calculated using the sessions described above, with the total time spent in the learning management system (LMS; for the data in this study this was Blackboard) that is used by most students to access WeBWorK. The resulting graph comparing these two times is shown in Figure 2 (this is for a one term Calculus 1 course with 1085 students, offered in 2016 at the institution where this study takes place). Figure 2a shows the data for all students in the course, and reveals a cluster of outliers whose LMS time is far lower than the time calculated directly from the WeBWorK log files. The reason for this discrepancy is that it is possible for students to access WeBWorK without going through the LMS; in this case the time spent in WeBWorK is not recorded as time spent in the LMS. Another source of error is that not all the time spent inside the LMS is necessarily time spent on WeBWorK; for the course in question, the LMS is also used to post students' grades and to host a discussion forum. Despite these sources of noise, we find that removing all students whose LMS time is less than half of their total WeBWorK time yields a graph with relatively strong correlation (Pearson's r = 0.77), as shown in Figure 2c. This graph is obtained by removing all the students below the dotted line in Figure 2a, and the blue line is the best fit line through the remaining data points in Figure 2c. The process of removing the students below the dotted line in Figure 2a can also be visualized using Figures 2b-2d, where we show the distribution of time (in WeBWorK or in the LMS) both before and after removing the outlier students. In Figure 2b, we see a bimodal distribution of times, which becomes unimodal in Figure 2d after the data points below the dotted line in Figure 2a are removed.

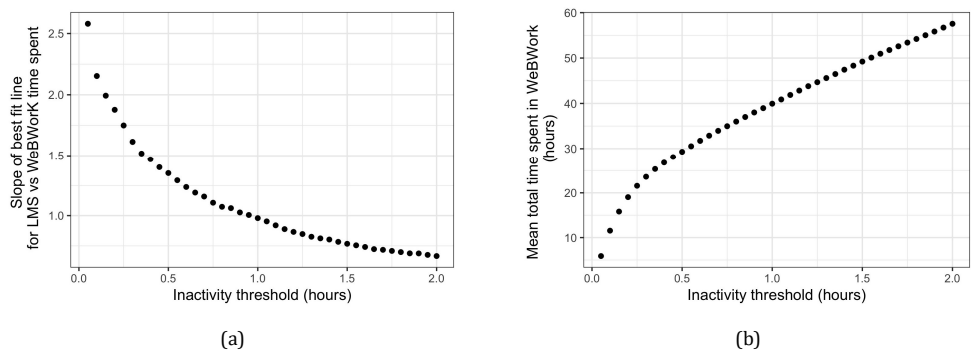

**Figure 3.** (a) The slope of the LMS time vs estimated WeBWorK time graph for different choices of the inactivity threshold. The threshold used in the remainder of the article is chosen so that the slope is closest to 1. (b) The average total estimated time spent in WeBWorK for different choices of the inactivity threshold. Based on the slope of the graph, each 0.1 hours change in the threshold corresponds to an approximately 5% change in the estimated total time spent in WeBWorK.

A source of uncertainty when calculating the total time spent inside WeBWorK is the choice of inactivity threshold. Here we have chosen the threshold of 0.95 hours, based on the criteria that this gives a slope closest to 1 when fitting a straight line to the data shown in Figure 2c. In Figure 3a, we show the values of the slope when fitting a straight line to the graph of LMS time vs WeBWorK time for different values of the inactivity threshold. Using the choice of 0.95 hours, we find that the blue line in Figure 2c has slope 1.005 (95% C.I. = [0.94,1.06]). Using a larger value for the inactivity threshold yields a larger estimate for the time a student spends in WeBWorK, thus decreasing the slope in Figure 2c. To illustrate the uncertainty in the estimate of total time in WeBWorK associated with the choice of inactivity threshold, in Figure 3b we show the average of total time in WeBWorK for all students in the course for choices of thresholds ranging from 0.1 to 2 hours. From the slope of the graph in Figure 3b we see that for every 0.5 hours increase in the threshold, the average total time in WeBWorK increases by approximately 10 hours (over the entire course). This corresponds to an error of approximately 25% in the average total time spent in WeBWorK for every 0.5 hours change in threshold, for the course chosen here and the choice of 0.95 hours for the threshold.

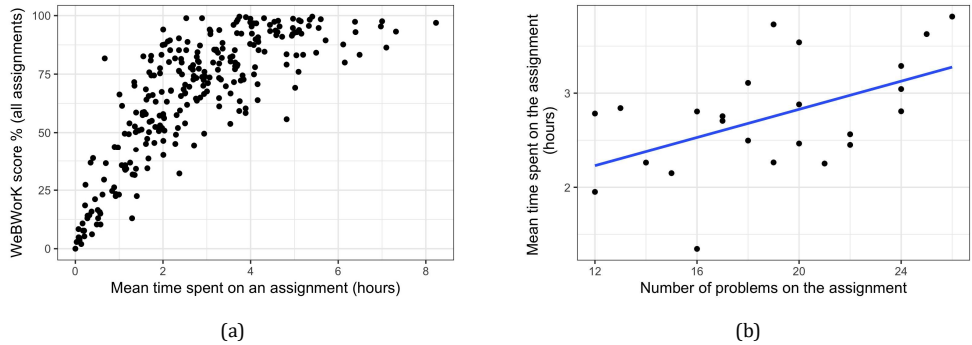

(a) (b)

**Figure 4.** (a) Total WeBWorK score % vs the estimated time spent in WeBWorK (average per assignment). (b) The estimated class average time spent on a WeBWorK assignment vs the number of questions on that assignment. As expected, there is a trend towards assignments with a larger number of questions requiring more student time on average (Pearson correlation = 0.52).

## 2.3. Validity of time estimates

Regardless of the choice of threshold, there will be uncertainty associated with the estimated time a student spends in WeBWorK. A larger inactivity threshold will count 'inactive' time between two sessions as active time, and thus overestimate the active time for a student that takes breaks between each completed WeBWorK problem that are long (e.g., 30-45 minutes) but still below the threshold. On the other hand, the procedure we outline above will underestimate the time spent on an assignment for a student who logs in and then spends more than 0.95 hours working on a single WeBWorK problem before making their first answer attempt on that problem.

Despite this uncertainty, we find that time estimated in this way is a good predictor of performance on those assignments, as shown in Figure 4a for a two-semester Calculus 1 course with 274 students in 2016-2017. The strong correlation between the average time spent on an assignment and overall WeBWorK performance helps to confirm that the method we propose can be used to indicate time that students spend on their assignments - at least relative to other students completing the same assignment even if there are uncertainties associated with the total time estimate. We also find evidence that our method captures the time students spend on an assignment (as a class average) relative to other assignments, as illustrated in Figure 4b, where we plot the time spent on an assignment (class average) against the number of problems on that assignment (same course as in Figure 4a; this course is used as the example for the remainder of this article).

As a further confirmation that our method correctly estimates the relative amount of time students spend on WeBWorK assignments, in Figure 5 we plot two other measurements that are expected to correlate with the average time spent on assignments. The first measurement is the answer to the survey question: 'On average, about how much time do you spend on a WeBWorK assignment for this course?' The survey question was administered to all 274 students in the class, with 96 respondents. The boxplots in Figure 5a illustrate the relationship between this self-reported time and the time estimated using the WeBWorK log files. We see that although the interquartile range of estimated times does not always fall within the self-reported range, the general trend is that (on average) students with a higher self-reported time spent also have a higher estimate of time spent from the log files. The graph in Figure 5b shows the average time spent on a WeBWorK assignment versus the average the number of answer attempts on an assignment (in these courses, multiple answer attempts per problem are allowed for almost all assignment problems). The number of attempts is strongly correlated with the estimate of time spent (Pearson's r = 0.81), with students who have a higher number of answer attempts at problems also spending a higher average time on assignments.

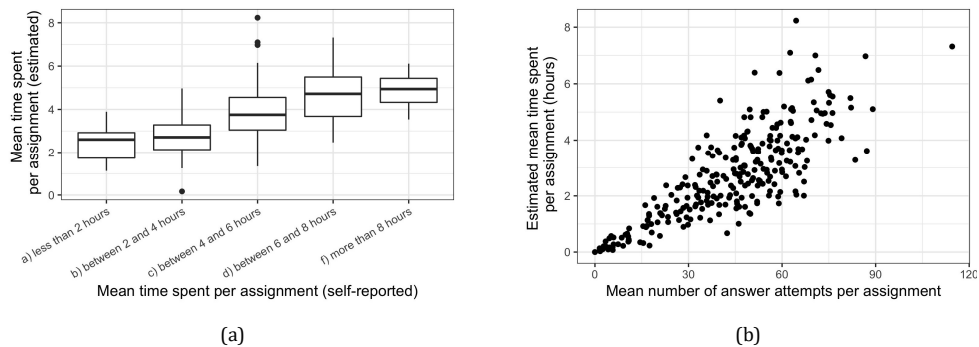

**Figure 5.** (a) Estimated WeBWorK time (average per assignment) vs self-reported time spent on WeBWorK assignments (answer to a multiple choice survey question). The solid lines are the medians and the boxes are the interquartile range. (b) The estimated mean WeBWorK time (per assignment) vs the average number of answer attempts made per assignment. Each point represents a student in the class. As expected, there is a strong correlation between the number of problem attempts made and the estimated time spent in WeBWorK.

## 3. Application: contrasting WeBWorK activity of groups of students

As an application of the analysis of WeBWorK log files and the estimate of time spent in WeBWorK discussed in the previous section, we now consider data from a single first year Calculus course and compare two groups of students within that course. The first group consists of those students with total WeBWorK assignment score less than 50%, and other group is the remainder of the class. Our motivation for considering these two groups of students is twofold. First, we use it as an opportunity to illustrate the type of information that can be extracted from WeBWorK log files, and how this information can be related to student performance on the WeBWorK assignments. Secondly, in the course considered here, there is a significant fraction of the class who struggle to complete the WeBWorK assignments, and this same group also has a high failure rate in the course; in 2016, 61% of students who had a total WeBWorK score of less than 50% over all assignments failed the course, compared to only a 6.3% failure rate in the remainder of the class. We are interested in knowing if it is possible to isolate some overall differences in behaviour that could explain the lower WeBWorK score of these students. These differences could then be used to address some of the factors that prevent students from successfully completing their WeBWorK assignments.

One metric which correlates with WeBWorK score is the total time spent on the assignment, as can be seen from Figure 4a. Given that students with lower WeBWorK scores spend less time on their assignments, a natural question which arises is whether the lower scoring students also score fewer points per hour than the higher scoring

**Table 1.** Comparison of two student groups on various metrics related to WeBWorK usage. On average, students in the lower scoring group (WW< 50%) attempt more points per hour and persist less on problems.

| Metric | WW [a] < 50% (n=64) Mean(SD) | WW ≥ 50% (n=209) Mean(SD) | Cohen's d |
|---|---|---|---|
| Points per hour | 7.2(4.7) | 5.1(2.2) | 0.7 |
| Number of problems attempted in one hour | 9.4(6.2) | 5.8(2.4) | 0.98 |
| Difficulty of problems attempted (scale: 0-100) | 35(2.1) | 34(0.34) | 0.29 |
| Time between first and last attempt on an incomplete[b] problem (hours) | 6.3(9) | 14(14) | -0.6 |
| Number of attempts on a problem attempted but not completed | 3.4(1.7) | 4.6(2.2) | -0.6 |

[a] WW - WeBWorK score % from total on all assignments. [b] an incomplete problem is one that has been attempted one or more times but never answered correctly

students. Surprisingly, we find that on average, the lower scoring students score more points per hour (Table 1). This same phenomenon can also be deduced from Figure 4a by noticing that the graph flattens out as one moves rightwards along the graph.

Given that students with low WeBWorK scores are capable of scoring more points per hour, a natural next question is whether these students are attempting more problems per hour when compared with the remainder of the class. The answer to this question is shown in Figure 6a-b, which shows that on average the lower scoring group attempt significantly more problems per hour. This suggests that the lower scoring students are employing a different strategy when it comes to scoring points on the assignments. One possible explanation for the observed trend in the data is that lower scoring students are selectively targeting problems which are easier, thus allowing them to attempt more problems per hour. To test this hypothesis, we give problems a difficulty rating based on the percentage of the class that attempted the problem but did not answer it correctly (scale: 0-100), and then calculate the average difficulty of the problems attempted by the low scoring students compared to the rest of the class. Surprisingly, we find that the lower scoring group attempt problems with a slightly higher average difficulty (Table 1). Because of this, we cannot conclude that students in the lower scoring group are able to score more points per hour by selectively targeting easier problems.

Another possible explanation for the observation that the lower scoring group attempt more problems per hour is that these students persist for less time when they are unable to correctly answer a problem. To test this hypothesis, we compare the low and high scoring groups on two measures of persistence. The first measure of persistence is the average time spent on a problem before giving up, which we estimate by using the average time between the first and last attempt on a problem that was attempted but not completed (attempted but not completed means that there was at least one answer submission made for that problem, but none of these answer submissions were fully correct). The second measure of persistence is the number of attempts at such a problem (that was attempted but not completed). For both of these measures of persistence the lower scoring group has a lower average measure of persistence (Table 1). The lower persistence of the lower scoring students, as measured in the average time spent on a problem before giving up, is shown in more detail in Figure 6c-d. This difference in persistence is consistent with the hypothesis that the lower scoring students are able to attempt more problems per hour by more hastily abandoning the problems which they are not able to complete.

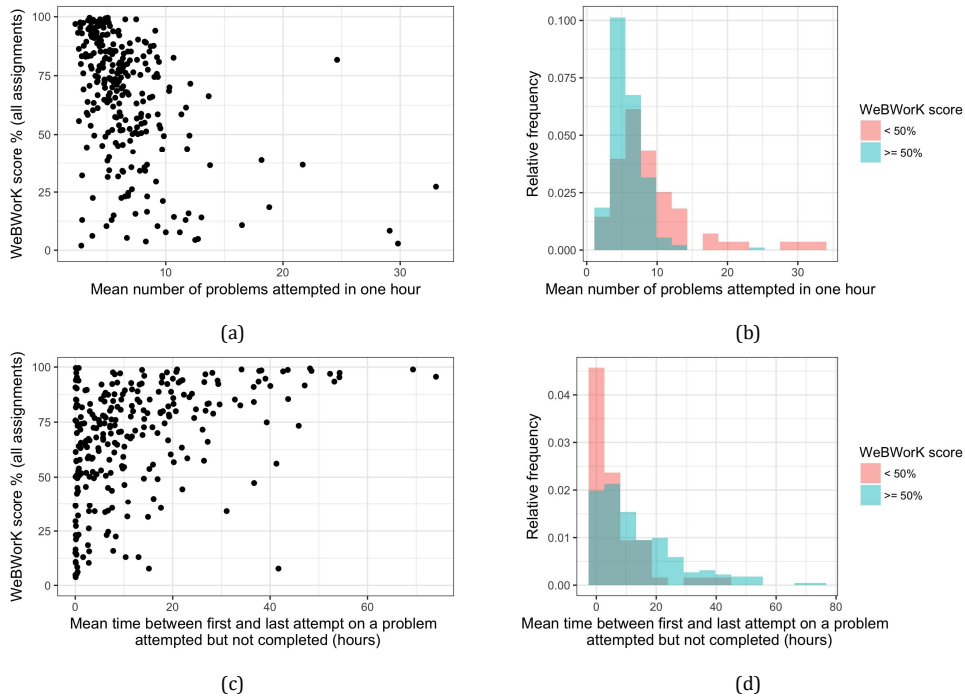

**Figure 6.** (a) WeBWorK score % vs average number of problems attempted per hour. Each point in the scatterplot represents a student in the class. The trend is that on average lower scoring students attempt more problems per hour. (b) Histogram contrasting the distributions of average number of problems attempted per hour for two groups of students: those with WeBWorK scores less than 50%, and those with WeBWorK scores greater than or equal to 50%. (c) WeBWorK score % vs the average time between the first and last attempt on a problem attempted but not completed (i.e. persistence time). The trend is that lower scoring students have shorter persistence times on average. (d) Distributions of persistence time for the high and low scoring groups of students.

**Table 2.** Comparison of the time spent on WeBWorK for high and low scoring student groups. On average, the lower scoring group spend less time on assignments. This is mainly due to spending spend fewer days on the assignment and have shorter activity sessions.

| Metric | WW [a] < 50% (n=64) Mean(SD) | WW ≥ 50% (n=209) Mean(SD) | Cohen's d |
|---|---|---|---|
| Time spent on an assignment (hours) | 0.91(0.72) | 3.3(1.4) | -1.89 |
| Session length (hours) | 0.37(0.15) | 0.48(0.12) | -0.89 |
| Number of login sessions per assignment | 2.3(1.5) | 6.9(2.7) | -1.89 |
| Time between first and last answer submission (days) | 1.1(0.91) | 2.2(1.3) | -0.84 |
| Time between sessions (hours) | 9.6(9.5) | 8.9(4.6) | 0.12 |
| First answer submission, in days before deadline | 2.3(1.4) | 3.2(1.8) | -0.55 |

[a] WW - WeBWorK score % from total on all assignments.

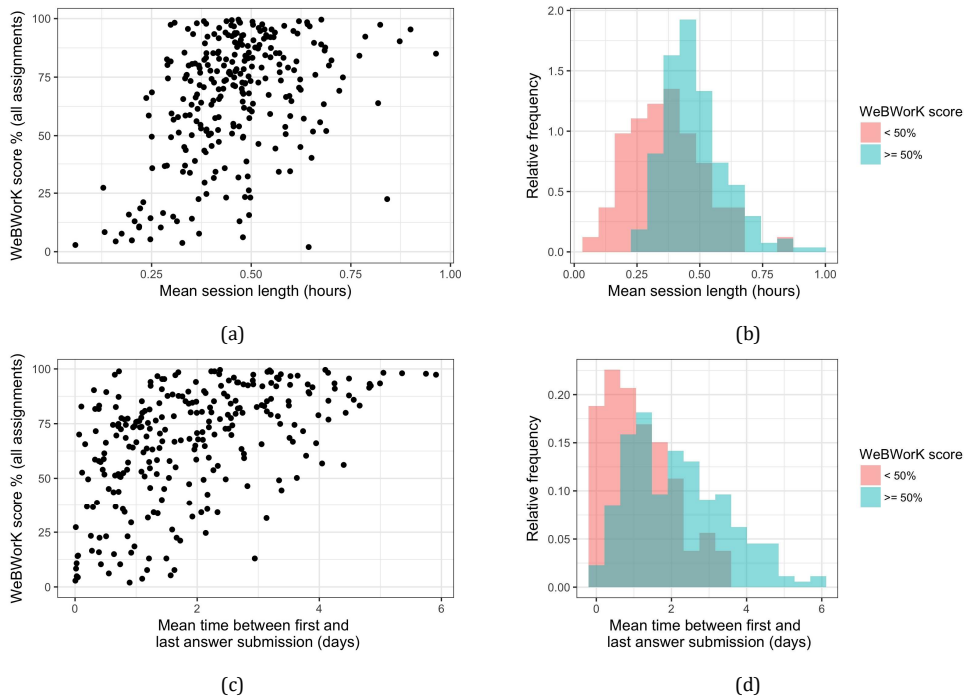

**Figure 7.** (a) WeBWorK score % vs average session length. The trend is that lower scoring students have shorter sessions on average. (b) Distributions of session length for students with WeBWorK scores less than 50% or WeBWorK scores greater than or equal to 50%. (c) WeBWorK score % vs average time between first and last answer submission (i.e. average start to finish time). The trend is that lower scoring students have shorter average start to finish time, and therefore concentrate their work over a shorter time period. (d) Distributions of start to finish time for the high and low scoring groups of students.

In addition to the questions involving total time on task and persistence considered above, a related question is how students use their activity sessions. For example, do lower scoring students have activity sessions that are shorter than other students, and do these activity sessions span a shorter period of time (when measuring the time between first and last session for a single assignment)? The answer to both of these questions is yes. The general trend is that lower scoring students have shorter average activity sessions (Table 2 and Figure 7a-b), and the average time between their first and last activity session for a single assignment is shorter (Table 2 and Figure 7c-d). Interestingly, we do not find any difference between the lower scoring group and the rest of the class in terms of the average number of hours between sessions (Table 2). Rather, it is the length of the sessions and the amount of time between the first and last session that accounts for the difference in time on task between the two groups. The shorter time between first and last session is likely to be related to lower scoring students starting their homework later relative to other students (last row of Table 2).

This difference in average starting times of approximately one day is most likely due to procrastination. The cause of the shorter activity sessions is less clear. It could be related to persistence, whereby students who encounter a sequence of problems that they are unable to complete decide to end their activity session out of frustration. Alternatively, it could be related to lower scoring students having a higher frequency of very short sessions that are caused by distractions that might be ignored by a more focused (higher scoring) student. Both of these possibilities would require deeper analysis of the log data, and both are worthy of further research.

## 4. Discussion

Our present work offers a scheme for learning from the WeBWorK log data that is routinely collected by the homework system, with a particular focus on what can be learned about the time that students spend interacting with the system. Our findings are compatible with existing results about online homework but have not previously been reported for WeBWorK. We have also attempted to collect studies from the various fields where WeBWorK studies can be found (mathematics education, practitioner focused publications, human computer interaction research, and the learning analytics community) to offer a comprehensive picture of the state of learning analytics approaches in WeBWorK.

By matching with independent lines of evidence (self-reports and other online activity), we have established that an estimate using log file data can characterize time spent on WeBWorK assignments. This information can lead to substantial further insight into student homework practices; to illustrate this, we have explored patterns in these time estimates to see that they align with previous results in online homework, where lower homework scores are associated with behaviours like procrastination, lower persistence, and more answer-seeking.

In future publications, we plan to further explore how WeBWorK log data can be used to understand student homework behaviours, and the implications of WeBWorK log data for teaching practices. One behaviour of particular interest is the case of students that scores very high on online homework but very low on exams. By combining the time estimates introduced here with patterns of answer submissions, it may be possible to show that some students are compromising their learning by taking homework shortcuts. Another potential application of WeBWorK log data would be to analyze the interaction between students' online forum discussions and their WeBWorK answer submissions. For example, by comparing the timing of answer submission with the timing of online discussions, it may be possible to observe students correcting a previously incorrect answer submission after seeking help on the discussion forum. As a third application of WeBWorK data, we plan to use the points per hour and other measures of performance on early WeBWorK assignments to identify students at risk of failing the course.

The other major area of application of WeBWorK log data is to improve the implementation of online assignments. For example, the relative class time spent on each assignment can be used to deduce the relative difficulty and workload at the scale of assignments (whereas the system currently only exposes instructors to difficulty measures on a per-question basis) and suggest adjustments for subsequent offerings of the course. In addition to helping with reflections on relative assignment difficulty, this type of information can be valuable to instructors as they plan assignments and how to support study skills in their courses. As has been suggested elsewhere (e.g., [3]), it is likely that additional tools like dashboards will be needed to make this information more easily available to instructors. Our work demonstrates that this information can be successfully built by drawing on existing sources of data captured by WeBWorK. Any further tools developed would stand to benefit instructors in a system that many already use rather than require a new homework system to participate. This type of scaling for an existing system has been relatively rare among learning analytics tools [19], but shows great promise in this case due to the maturity and widespread use of WeBWorK.

## Acknowledgements

This work was supported by the Carl Wieman Science Education Initiative at the University of British Columbia. The authors would like to thank Costanza Piccolo for feedback on the manuscript.